\documentstyle[amstex,12pt,a4,twoside]{amsart}
\def\myauthor{Ph\`ung H{{\accent"5E o}\kern-.28em\raise.2ex\hbox{\char'22}\kern-.20em}
 H{a\kern-.370em\raise.16ex\hbox{\char'47}\kern.1em}i}
\newcommand{\BBb}[1]{{\Bbb #1}}

\def\ot{\otimes}

\newcommand{\bbas}{\begin{eqnarray*}}
\newcommand{\eeas}{\end{eqnarray*}}
\newcommand{\bbar}{\begin{array}}
\newcommand{\eear}{\end{array}}
\newcommand{\bbs}{\begin{displaymath}}
\newcommand{\ees}{\end{displaymath}}
\newcommand{\bb}{\begin{equation}}
\def\ee{\end{equation}}
\def\eea{\end{eqnarray}}

\def\bba{\begin{eqnarray}}
\newtheorem{thm}{Theorem}[section]
\newtheorem{lem}[thm]{Lemma}

\newtheorem{rem}[thm]{Remark}
\newtheorem{edl}[thm]{Theorem}

\newtheorem{cor}[thm]{Corollary}

\def\det{\mbox{\rm det{}}}

\def\Im{\mbox{\rm Im}{}}

\def\End{\mbox{\rm End{}}}

\def\H{{\cal H}}

\def\eee{\rule{.75ex}{1.5ex}\vskip1ex}

\def\proof{{\it Proof.\ }}
\renewcommand{\dim}{\mbox{\rm dim}}

\def\rref#1{(\ref{#1})}

\font\Fraktur=eufm10 scaled\magstep1          % for display- and textstyle
   \newcommand{\fraktur}[1]{\mbox{\Fraktur #1}}  %
   \font\Fraktu=eufm7 scaled\magstep1            % for scriptstyle
   \newcommand{\fraktu}[1]{\mbox{\Fraktu #1}}    %
   \font\Frakt=eufm5 scaled\magstep1             % for scriptscriptstyle
  \newcommand{\frakt}[1]{\mbox{\Frakt #1}}      %
   \def\fr#1{\mathchoice{\fraktur {#1}}            % displaystyle
                        {\fraktur {#1}}            % textstyle
                        {\fraktu {#1}}             % scriptstyle
                        {\frakt {#1}}  }           % scriptscriptstyle

\newcommand{\Ss}{\fr S}

\newcommand{\Comod}{\mbox{Comod}}
\def\part{\vdash}

\def\part{\vdash}

\newcommand{\bfe}{{\mathbf e}}
\newcommand{\bfp}{{\mathbf p}}
\newcommand{\bfm}{{\mathbf m}}
\newcommand{\bfs}{{\mathbf s}}
\newcommand{\bfk}{{\mathbf k}}
\newcommand{\bfLambda}{\mbox{\boldmath$\Lambda$\unboldmath}}
\newcommand{\bflambda}{\mbox{\boldmath$\lambda$\unboldmath}}

\def\ictpthanks{ }
\def\myaddress{International Center for Theoretical Physics\\ P.O. Box 586, 34100 Trieste, Italy}
\def\mythanks{The work was done during the author's stay at the International Center for Theoretical Physics, Trieste, Italy. He would like to thank the Center for its hospitality and financial support.}
\def\myabstract{We prove that the dimension of the homogeneous components of the quadratic algebras associated to a Hecke operator is a P\'olya frequency (P-) sequence. This result enables us to give some characterizations on the Poincar\'e series of  these quadratic algebras.}

\def\amshead{
\title{ Poincar\'e series of quantum spaces associated to Hecke operators }
\author{\scshape\myauthor}
\address{\myaddress \footnote{On leave from Hanoi Institute of Math, Hanoi, Vietnam.}}
\thanks{\mythanks}
\keywords{Hecke operator, total positivity, Poincar\'e series.}
\email{phung@@ictp.trieste.it}
\dedicatory{\hfill \normalsize\Fraktur Herrn Prof. Dr. Bodo Pareigis zum 60. Geburtstag gewidmet}
\begin{abstract}\myabstract\end{abstract}
\subjclass{Primary 16W30,17B37 , Secondary 17A45, 17A70}
\maketitle
}
\begin{document}
\amshead
%\ictphead
\bibliographystyle{plain}
%\begin{flushright}
%\hfill \normalsize\sl Herrn Prof. Dr. Bodo Pareigis zum 60. Geburtstag gewidmet
%\end{flushright}
\section{Introduction}
Subject of this work are the quadratic algebras associated to a Hecke operator. Hecke operator on a vector space $V$ is a linear operator on $V\ot V$ that satisfies the Yang-Baxter equation and the Hecke equation $(x+1)(x-q)=0$ and hence induces a representation of the Hecke algebra of the general linear group on the tensor power $V^{\ot n}$ of $V$ for every $n\geq 2$. Given a Hecke operator $R$ on a vector space $V$, one can define the quantum anti-symmetric and symmetric tensor algebras on $V$, which, in case $R$ is the usual permuting operator, are the usual anti-symmetric (i.e. exterior) and symmetric tensor algebras over $V$. Thus one can think of these quadratic algebras as the  anti-symmetric and symmetric tensor algebras on certain quantum space. Then one can define an algebra that coacts universally upon this pair of quadratic algebras, which turns out to be a quadratic algebra, too. This algebra can be considered as the function algebra over the matrix quantum semi-group of ``symmetries'' of the quantum space mentioned above\cite{manin1,lyu,gur1,sud}.

In this work we study the Poincar\'e series of the above quadratic algebras. The first attempt in this study was made by Lyubashenko, who proved that (in case of $q=1$) if $P_\Lambda (t)$ is a polynomial then it is a reciprocal polynomial iff $R$ is closed \cite{lyu}. This result was generalized by Gurevich for arbitrary $q$ \cite{gur1}. He also gave an example of an even Hecke operator whose rank differs from  the dimension of the vector space it is defined on. Here we show that the rank should be less than the dimension of the corresponding vector space. We also study the case of odd-even Hecke operator and analogously show that that its super rank should have the weight (the sum of two components of the rank) less than the dimension of the vector space it is defined on. Our main theorem (Theorem \ref{dl35}) states that the Poincar\'e series of a quantum space associated to a Hecke operator is a rational function having negative roots and positive poles.
 
Using the theory of symmetric functions we give a combinatorial formula for the dimension of irreducible representation of the matrix quantum semi-group as functions on the dimension of homogeneous components of the quantum space. Whence we recover the relation between the Poincar\'e series of the quantum space and of its matrix quantum semi-group, which was found in \cite{ph97}.

In this work, a theorem of Edrei on P\'olya frequency (or P-) sequences \cite{edrei} plays an important role, from which almost all results follow. Our method bases on the theory of Schur symmetric functions \cite{mcd1,mcd2} as well as the quantum version of Schur's double centralizer theorem \cite{ph97}.

\section{Quantum spaces and quantum semi-groups.}
Let $\BBb{K}$ be a field of characteristic zero, $V$ be a
 vector space over $\BBb{K}$ of finite dimensional $n$. A Hecke operator on $V$ is a $\BBb{K}$-linear operator $R : V \otimes V \longrightarrow V \otimes V$, satisfying the following relations:
\bba\label{eq1} &
(R \otimes I)(I \otimes R)(R \otimes I) = (I\otimes R)(R \otimes I)(I \otimes R),& \\ \label{eq2} &
(R + I)(R - q) = 0,  0 \neq q \in \BBb{K}.&
\eea
The name Hecke operator comes from the fact, that an operator, satisfying the above equations, induces a representation of the Hecke algebra $\H_n$ on $V^{\otimes n}, \forall n \geq 2.$
The reader is referred to \cite{dj1} for the definition and main properties of Hecke algebras.
 
A quadratic algebra is a pair $(A,V)$ consisting of a finite dimensional vector space $V$ and a factor algebra $A$ of the tensor algebra $T(V)$ by an ideal, generated by elements of $V \otimes V$. Given a Hecke operator $R$ on $V \otimes V$, one can define the following quadratic algebras.
\bba\label{eq3}&
(\Lambda, V) = (\Lambda_R, V):= (T(V)/<\Im(R + 1) >, V)&\\
\label{eq4}&
(S, V) = (S_R, V): = (T(V)/<\Im (R-q)>, V)&\\
\label{eq5}&
(E, V^* \otimes V) = (E_R, V^* \otimes V) := (T(V^* \otimes V)/<\Im(\overline{R}-1)>,V^* \otimes V), &
\eea
where $ \overline{R} = \theta (R^{*-1} \otimes R) \theta,$ acting on $ V^* \otimes V \otimes V^* \otimes V,\ \theta $ denotes the operator, that interchanges the  elements in the $ 2^{\rm nd} $ and $ 3^{\rm rd} $ positions of the tensor products.
Fix $V$ and $R,$ we shall refer to these algebras as $\Lambda, S$ and $E$, respectively. $\Lambda$ and $S$ are considered as the ``anti-symmetric and symmetric tensor algebras'' upon certain quantum space, while $E$ is considered as the algebra of functions on the ``semi-group of symmetries of this quantum space'', i.e., the matrix quantum semi-group.

$\Lambda, S$ and $E$ are graded algebras, let $\lambda_i, s_i, e_i$ be the $\BBb{K}$-dimension of their $i^{\rm th}$ homogeneous components, respectively. For a quadratic algebra $(A,V)$, its Poincar\'e series is defined to be $P_A (t) = \sum_{i=0}^\infty (\dim_{\BBb{K}} A_i)t^i$, where $A_i$ is the $i^{\rm th}$-homogeneous component of $A$. Thus, for example, $P_{\Lambda} (t) = \sum_{i=0}^\infty \lambda_it^i$.
The Poincar\'e series of $\Lambda$ and $S$ satisfies the following relation \cite{gur1}:
\begin{equation}\label{eq6}
P_\Lambda (t) P_S (-t) = 1.
\end{equation}

$E$ is a bialgebra and $\Lambda, S$ are (right) $E$-comodules. In fact $E$ can be defined as an algebra that universally coacts on $\Lambda$ and $S$ \cite{manin1}. It then immediately implies that $E$ is a bialgebra. Also, $V$ is an $E$-comodule and hence $V^{\otimes n}$ is an $E$-comodule. The coproduct on $E$ satisfies $\Delta(E_n) \subset E_n \otimes E_n$, hence $E_n$ is a coalgebra. We also have $\delta (V^{\otimes n}) \subset V^{\otimes n} \otimes E_n,$ i.e., $V^{\otimes n}$ is an $E_n$-comodule. Therefore $E_n^*$ is an algebra, that acts on $V^{\otimes n}$ (from the left). 

On the other hand, as has been mentioned, $R$ induces a representation, say $\rho_n$, of the Hecke algebra $\H_n$ on $V^{\otimes n}$. We shall write this action from the right. It is now time to state the fact, that plays a crucial role in the study of $E$.

\begin{edl}\label{dl11}{\rm \cite{ph97}}
The actions of $E_n^*$ and $\H_n$ on $V^{\otimes n}$ are centralizers of each other in $\End_{\BBb{K}} (V^{\otimes n})$. In other words, the following isomorphisms hold
\bba\label{eq7}&
E_n^* \cong \End_{\H_n} V^{\otimes n},&\\
\label{eq8}&
\rho_n (\H_n) \cong \End_{E_n^*} V^{\otimes n}.&
\eea\end{edl}
Since $ \H_n$ is semi-simple and hence isomorphic to the direct product of matrix rings over $\Bbb K$ (by comparing dimensions), these equations  imply the following facts.
\begin{itemize}
\item[i)] $E_n^*$ is semi-simple, hence $E_n$ is cosemi-simple.
\item[ii)] Simple $E_n^*$-modules are of the form $\Im\rho_n (e)$, where $e$ is a primitive idempotent of $ \H_n$. Two primitive idempotents define isomorphic modules iff they belong to the same minimal two-sided ideal.
\item[iii)] \begin{equation}\label{eq9}
E_n^* \cong \oplus\End_{A_i} V^{\otimes n} A_i,
\end{equation}
where $A_i$ are minimal two-sided ideals of $ \H_n$.
\item[iv)] $\End_{A_i}(V^{\otimes n} A_i)$ is a matrix ring of dimension equal to the square of the $\BBb{K}$-dimension of $\Im\rho_n (e)$ for any primitive idempotent $e \in A_i$.
\end{itemize}

On the other hand, it is known that minimal two-sided ideals of $ \H_n$ can be indexed by partitions of $n$ \cite{dj2}. Let $\lambda$ be a partition of $n$ (in notation $\lambda \vdash n$), that means $\lambda = (\lambda_1, \lambda_2,\ldots), \lambda_1\geq \lambda_2\geq\ldots ,\sum\lambda_i = n$, and $A_\lambda$ be the corresponding two-sided minimal ideal, then let $M_\lambda$ be the simple module of $E_n^*$ defined by any primitive idempotent from $A_\lambda$. Let $m_\lambda$ denote the dimension of $M_\lambda$ over $\BBb{K}$. Then we have, according to (iv) above,
\begin{equation}\label{eq10}
e_n = \dim E_n = \dim E_n^* = \sum_{\lambda\vdash n} m_\lambda^2.
\end{equation}
Notice that by definition $m_\lambda\geq 0$ for all $\lambda$. This remark together with Theorem \ref{dl22} implies almost all the results of this paper. 

Since $R$ satisfies the Yang-Baxter equation, it induces a coquasitriangular structure on $E$ and, therefore, a braiding in the category $E$-\Comod\ of $E$-comodules. In particular, the braiding on $V\otimes V$ is precisely $R$. The reader is referred to \cite{js1,js2} for the general theory of braided categories and to \cite{l-t} for the particular case of $E-$\Comod .

\section{Symmetric functions and P\'olya frequency sequences}

In this section we recall some properties of symmetric functions needed in our context. The reader is referred to the book of Macdonald \cite{mcd1,mcd2} for an introduction into the theory of symmetric function. We then recall a theorem, due to Edrei \cite{edrei}, that plays  the main role in our study.

We study the symmetric functions in a countable set of variables $x_1, x_2,\ldots$. We shall closely follow Macdonald's book \cite{mcd1,mcd2}, up to some changes of notation for our convenience.
% make some changes of notation, comparing with Macdonald's ones \cite{mcd1,mcd2}, for our convenience. 
The monomial symmetric functions will be denoted by $\bfk_\lambda, \lambda\in{\cal P}$, where ${\cal P}$ denotes the set of partitions. The elementary, complete and Schur functions will be denoted by $\bflambda_\lambda , \bfs_\lambda$ and $\bfm_\lambda , \lambda \in {\cal P}$, respectively. Recall that $\bflambda_\lambda = \bflambda_{\lambda_1}\bflambda_{\lambda_2}\cdots\bflambda_{\lambda_n} , \bfs_\lambda = \bfs_{\lambda_1} \bfs_{\lambda_2}\cdots \bfs_{\lambda_n}$, for $\lambda = (\lambda_1, \lambda_2, \cdots\lambda_n)$, where $\bflambda_r$ (resp. $\bfs_r$) is the coefficient of $t^r$ in $\prod\nolimits^\infty_{i=1} (1 + x_i t)$ (resp. $\prod\nolimits^\infty_{i=1} (1 - x_i t)^{-1}$). $\bfm_\lambda$ is given by the equality (we make use of the convention that $\bfs_r = \bflambda_r =0$ whenever $r<0$):
\begin{equation}\label{eq11}
\bfm_\lambda = \det\parallel\bfs_{\lambda_i - i +j}\parallel_{1\leq i, j\leq l(\lambda)} = \det\parallel\bflambda_{\lambda_i -i +j}\parallel_{1\leq i, j\leq l(\lambda ')}\end{equation}
 where $l(\lambda ')$ denotes the length of $\lambda$, i.e., the cardinal of nonzero parts of $\lambda,$ $\lambda '$ denotes the conjugate partition to $\lambda$. The set $\bfk_\lambda$, (resp. $\bflambda_\lambda , \bfs_\lambda , \bfm_\lambda), \lambda \in {\cal P}$ forms a basis for $\bfLambda$ -- the rings of symmetric functions ([loc.cit.] I.2,3). In particular, we have $\bfm_\lambda=\sum_\mu K_{\lambda}^\mu\bfk_\mu$, the coefficients $K_{\lambda}^\mu$, called Kostka number, are nonnegative integers ([loc.cit.] I.4,5).

We also make use of skew Schur functions $\bfm_{\lambda/\mu}$ defined as follows:
\begin{equation}\label{eq12}
\bfm_{\lambda/\mu} : = \det\parallel \bfs_{{\lambda_i} - {\mu_i} -i+j}\parallel_{1\leq i,j\leq l(\lambda)} = \det\parallel \bfe_{{\lambda_i '} - {\mu_j '} -i+j}\parallel_{1\leq i,j\leq l(\lambda')}.\end{equation}
For $\mu = 0$, it becomes $\bfm_\lambda$.
$\bfm_{\lambda / \mu}$ can be expressed in terms of $\bfm_\lambda$. Let $c_{\lambda\mu}^\gamma$ be the (integer) number such that $\bfm_{\lambda/\mu} = \sum c_{\lambda\mu}^\gamma \bfm_\gamma$, then
\begin{equation}\label{eq13}
\bfm_{\lambda /\mu} = \sum c_{\mu\gamma}^\lambda \bfm_\gamma.
\end{equation}
$c_{\mu\gamma}^\lambda =0$ unless $|\lambda | = |\mu | + |\gamma |$ and $\lambda\supset\mu$ (i.e., $\lambda_i\geq\mu_i , \forall i$).
Moreover $c_{\mu\gamma}^\lambda$ are nonnegative ([loc.cit.] I.5).\\[1ex]

In the previous section, we have seen that each partition $\lambda$ defines a simple $E$-comodule $M_\lambda$ (which can be zero). If $M_\lambda$ is zero, we can assume it appears in the decomposition of an $E$-comodule with an arbitrary multiplicity. In particular, consider the multiplicity of $M_\gamma$ appearing in the decomposition of $M_\lambda\ot M_\mu$, if one of $\lambda,\mu,\gamma,$ equals zero, then the multiplicity can be chosen arbitrary. We show that this multiplicity can be chosen so that it depends only on $\lambda,\mu,$ and $\gamma.$
\begin{lem}\label{lem21}Let $M_\lambda, \lambda\in {\cal P}$ be the simple $E$-comodules, some of them may be zero.
Let $M_\lambda \otimes M_\mu = \bigoplus d_{\lambda\mu}^\gamma M_\gamma$. Then $d_{\lambda\mu}^\gamma$ can be chosen to be equal to $c_{\lambda\mu}^\gamma$.
\end{lem}
\proof
 Let $e_{\lambda}$ (resp. $e_\mu)$ be a primitive idempotent of $\H_m$ (resp. $\H_n)$ where $m = |\lambda |$ (resp. $ n = |\mu |)$. Embed $\H_m$ and $\H_n$ into $\H_{m+n}$ in such a way that the generators $T_i, 1\leq i\leq m,$ of $\H_m$ are mapped to the generators $T_i$ of $\H_{m+n}$ and the generators $T_j, 1\leq j \leq n$, of $\H_n$ are mapped to the generators $T_{j+m}$ of $\H_{m+n}$. Then $d_{\lambda\mu}^\gamma$ can be chosen to be equal to the cardinal of primitive idempotents belonging to the minimal two-sided ideal indexed by $\lambda$ of $H_{m+n}$, that appear in the decomposition of $e_\lambda e_\mu$ into sum of primitive idempotents. Thus $d_{\lambda\mu}^\gamma$ is an integer and does not depend on $q$, hence we can set $q=1$. In this case $\H_n = \BBb{K} [\Ss_n]$, hence $d_{\lambda\mu}^\gamma = c_{\lambda\mu}^\gamma$ (cf. \cite{green1}, Chapter 3).\eee

Since $\{ \bflambda_r | r\geq 0\}$ are algebraically independent, we can assign to them any set of values from $\BBb{K}$ to obtain a ring homomorphism $\bfLambda \longrightarrow\BBb{K}$. In particular, we can set $\bflambda_r = \lambda_r$ -- the dimension of $\Lambda_r$ defined in the previous section. By virtue of equation \rref{eq6}, $\bfs_r = s_r = \dim  S_r$. Furthermore, since the dimension of the tensor-product of two modules is equal to the product of their dimension, by Lemma \ref{lem21}, $\bfm_\lambda = m_\lambda = \dim M_\lambda$.
Thus, by the definition of $\bfm_\lambda$, we have
\begin{equation}\label{eq14}
m_\lambda = \det |s_{\lambda_i-j+i}|_{1\leq i,j\leq l(\lambda)} = \det|\lambda_{\lambda'_i-j+i}|_{1\leq i,j\leq l(\lambda ')}.
\end{equation}

Since $m_\lambda$ is the dimension of a vector space, it is nonnegative.
We make use of the following result from complex analysis to get information on $\lambda_r$.
\begin{edl}\label{dl22}{\rm \cite{edrei}}
Let $a(z) = \sum_{i=0}^\infty a_i z^i $ be a formal series of complex coefficient, $a_0=1$. Then the infinite matrix $\parallel a_{i-j} \parallel _{0\leq i,j}$ (we make the convention that $a_i =0$ whenever $i<0$) is totally positive, that is, all its minor of any degree $r\geq 1$ is nonnegative, if and only if $a(z)$ is of the form
\begin{equation}\label{eq15}
a(z ) = e^{ \gamma z} \frac{\prod_1^{\infty} (1+t_i z )}{\prod_1^{\infty} (1-u_i z)},
\end{equation}
with $\gamma\geq 0, t_i \geq 0, u_i\geq 0$ and $\sum (t_i + u_i) \leq \infty$. This equation is understood to be an equation of  two convergent complex functions in a neighborhood of zero.
\end{edl}
 
The series $ a(t) = \sum_{i=0}^\infty a_i t_i$ {\it generates} the sequence $a_0, a_1, \ldots, $ which is called a P\'olya frequency or (P-) sequence if $a(z)$ satisfies the condition of Theorem \ref{dl22}. 

\begin{lem}\label{lem23}
Let $R$ be a Hecke operator. Then $\lambda_0, \lambda_1, \ldots,$ is a P-sequence.
\end{lem}
\proof
Indeed, any minor of $\parallel \lambda_{i-j} \parallel _{0\leq i,j}$ has the form: $\det \parallel \lambda_{{\mu_i}-{\nu_j}} \parallel_{ 1\leq i, j\leq r}$ where $\mu , \nu$ are two sequences of nonnegative distinct integers $r=l(\mu )$. We have, according to \rref{eq12},
\begin{equation}\label{eq16}
\det \parallel \lambda_{{\mu_i} - {\nu_j}} \parallel_{1\leq i, j \leq r} = \bfm_{{\mu + \delta}/{\nu +\delta}} = : m_{{\mu +\delta}/{\nu + \delta}},
\end{equation}
where $\delta = (r-1, r-2, \cdots, 0)$. Since $m_{\mu / \nu}$ is a linear  combination of $m_\gamma$ with nonnegative coefficients, it is itself nonnegative.\eee

A P-sequence is called PP-sequence if all the minors $\parallel\lambda_{{\mu_i} - {\nu_j}}\parallel_{1\leq i, j \leq r}, \nu \subset\mu,$ are strictly positive. We have
\begin{lem}\label{lem34} Let $a(z)$ be of the form in \rref{eq15}. If any of the following conditions is fulfilled: (1) $\gamma> 0$, (2) $a_i>0, \forall i$, (3) $b_i>0, \forall i$, then $a(z)$ is a generating function for a PP-sequence.\end{lem}
\proof The Lemma follows from the following observations.

1. Let $x,y$ be two countable sets of positive reals, then $\bfk_\lambda(x,y)\geq \bfk_\lambda(x)$,  $\bfk_\lambda(x,y)\geq \bfk_\lambda(y)$.

2. Analogously we have $\bfm_\lambda(x,y)\geq \bfm_\lambda(x)$,  $\bfm_\lambda(x,y)\geq \bfm_\lambda(y)$ (\cite{mcd1,mcd2}, (5.9)).

3. From the first observation, we see that, for $a(z)=\prod_1^\infty(1+a_iz), a_i>0$, $\bfk_\lambda>0$ for all $\lambda$, hence $\bfm_\lambda > 0, \forall \lambda$. Therefore, if $a(z)=\prod_1^\infty(1-a_iz)^{-1}, a_i>0$ we also have $\bfm_\lambda > 0$ (using the involution $\omega$ [loc.cit.], I.3).

4. From the second observation, we see that if $a(z)$ generates a P-sequence and $b(z)$ generates a PP-sequence then $c=a\cdot b$ generates a PP-sequence.

5. $e^{\gamma z}$ generates a PP-sequence (cf. [loc.cit], I.3, Ex. 5).

3.,4.,5. imply the assertion of the Lemma.\eee
As a direct consequence of this lemma we have
\begin{edl}\label{dl35} Let $R$ be a Hecke operator then $P_\lambda(t)$ is a rational function.\end{edl}
\proof It is sufficient to show that there exists a partition $\lambda$, such that $M_\lambda$ is zero. Were $M_\lambda\neq 0$ for all $\lambda$, then $\rho_n$ in \rref{eq8} were injective for all $n$, by comparing the dimension we would have a contradiction.\eee
%\begin{rem}\label{rem34}\rm Concluding this section we would like to recall a fact, which will be of use in the next sections, that Schur functions are linear combinations of monomial functions with positive integer coefficients \cite{mcd1}, hence Schur functions defined on a set of nonnegative number are nonnegative.\end{rem}

\section{ The Quasi-even Case}
A Hecke operator $R$ is called quasi-even if the associated quadratic algebra $\Lambda=\Lambda_R$ is finite dimensional, i.e., if $P_\Lambda(t)$ is a polynomial. $R$  is called even Hecke operator (or symmetry)  if the dimension of the highest homogeneous component of $\Lambda$ is equal to 1. The degree of $P_\Lambda$ is called the rank of $R$.

Let now $R$ be a quasi-even Hecke operator. Then, according to Theorem \ref{dl35}, $P_\Lambda(t)$ has only non-positive roots %(cf. \cite{brenti} Theorem 4.5.3), and since $P_\Lambda(0)=1$, 
\begin{equation}\label{eq17} P_\Lambda(t)=\prod_{i=1}^r(1+t_it), 0\neq t_i\in\BBb{R}.\end{equation}
Hence, according to \rref{eq11}, $m_\lambda\neq 0$ if and only if $l(\lambda)>r$. Thus, simple $E$-comodules can be indexed by $\{ \lambda\in {\cal P}|l(\lambda)\leq r\}$.
\begin{edl}\label{dl31} Let $R=R_q$, $R'=R'_{q'}$ be quasi-even Hecke symmetries over $V$ and $V'$, of ranks $r$ and $r'$, respectively. Let $E$ and $E'$ be the associated  bialgebras. Then the categories $E$-\Comod\ and $E'$-\Comod\ are equivalent as braided abelian categories if and only if $q=q'$ and $r=r'$.\end{edl}
\proof
Let ${\cal F}$ be the equivalence, then ${\cal F}(V)=V'$, hence ${\cal F}(V^{\ot n})=V^{'\ot n}$. According to \rref{eq8}, 
\bbs\End^E(V^{\ot n})\cong\End_{E^*_n}(V^{\ot n})=\bigoplus_{\lambda\part n,l(\lambda)\leq r}A_\lambda,\ees
 where $A_\lambda$ denotes the minimal two-sided ideal of $\H_n$, associated to $\lambda$. Therefore $q=q'$ and $r=r'$.\eee

Remark that the dimension of $V$ plays no role here.

Assume now that $R$ is an even Hecke operator, i.e., $\lambda_r=1$. In this case one can show that $P_\Lambda(t)$ is a reciprocal polynomial, i.e., $\lambda_k=\lambda_{r-k}, k=0,1,\ldots,r$, \cite{gur1}. Since $\lambda_r=1, \prod_1^r t_i=1$, hence $\lambda_1=\sum_it_i\geq r$. Notice that $\lambda_1=\dim V$. Thus, we have
\begin{edl}\label{dl32} Let $R$ be an even Hecke operator of rank $r$ over the vector space $V$ of dimension at least 2. Then $2\leq r\leq n$, (cf. \cite{gur1}, Corollary 4.5).\end{edl}
Notice that $r=1$ means $R=qI$, which is not even unless $\dim V=1$. If $r=\dim(V)$, then $t_i=1, \forall i$. Thus we have
\begin{cor}\label{cor33}Let $R$ be an even Hecke operator over $V$ such that the rank of $R$ is equal to the dimension of $V$ then $P_\Lambda(t)=(1+t)^r$, hence $P_S(t)=(1-t)^{-n}$ and by Equation \rref{eq23} $P_E(t)=(1-t)^{-n^2}$.\end{cor}

\begin{rem}\rm As we have already seen, examples of strictly quasi-even Hecke operators serve the identity operators. On the other hand, it is proved in \cite{ph97b} (and can also be derived from Lemma \ref{lem21}), that $\Lambda _r^{\ot 2}$, where $r$ is the rank of $R$, is a simple $E$-comodules, hence the braiding on  $\Lambda_r^{\ot 2}$ is the identity operator. For more examples of strictly quasi-even Hecke operators one can take the Hecke sum (see Section \ref{secE}) of an identity operator and an even Hecke operator. I conjecture that quasi-even Hecke operators appear in this form.

For an even Hecke operator, we have shown that $P_\Lambda$ should be a reciprocal polynomial having only negative roots. I conjecture that any reciprocal polynomial having negative roots and leading coefficient equal to 1 should have the form $P_\Lambda(t)$ for some even Hecke operator. This fact, for the case of polynomials of degree 2 and 3, follows from a result due to Gurevich \cite{gur1}.\end{rem}

\section{The quasi-odd-even case}
If $P_\Lambda (t)$ is not a polynomial $R$ is called quasi-odd-even Hecke operator. We have
%a rational function. According to \rref{eq15}(cf. \cite{brenti} Theorem 4.5.3),
\begin{equation}\label{eq18}
P_{\Lambda} (t) = \frac{\prod_1^m (1+ t_i t)}{\prod_1^n (1-u_j t)},
\end{equation}
where $t_i>0, u_j>0$. $(m,n)$  is called the super rank of $R$. $R$ is called odd-even Hecke operator if $\prod t_i = \prod u_i =1$. Not much is known about (quasi-) odd-even Hecke operators. 
%On the other hand, it is not known if any Hecke operator is quasi-odd-even.

For expressing dimension of simple $E$-comodules we use the super Schur (or Hook) functions. Let $x=(x_1, x_2, \cdots ,x_m)$ and $y=(y_1, y_2, \cdots , y_n)$ be two sets of variables. The super Schur function $\bfm_\lambda (x/y)$ is defined as follows:
\begin{equation}\label{eq19}
\bfm_\lambda (x/y):= \sum_{\mu\subset\lambda} \bfs_\mu (x_1, x_2, \cdots, x_m) \bfs_{\lambda '/\mu '} (y_1, y_2, \cdots ,y_n)
\end{equation}
This function was introduced by Berele and Regev \cite{br1}.  Then we have \cite{mcd2}
%\begin{equation}\label{eq20} 
$m_\lambda = \bfm_\lambda (t/u).$
%\end{equation}
Hence $m_\lambda \neq 0$ if and only if $\lambda\in\Gamma_{m,n}$ where $\Gamma_{m,n}$ denotes the set of hook-partitions:
\begin{equation}\label{eq21}
\lambda\in\Gamma_{m,n} \Longleftrightarrow \lambda_j \leq n \mbox{ for } j\geq m+1.
\end{equation}
Thus we have 
\begin{edl}\label{dl41}
Let $R$ be a quasi-odd-even Hecke operator then simple $E$-comodules can be indexed by $\lambda\in\Gamma_{m,n}$. The categories of comodules over the bialgebras $E$ and $E'$, associated with quasi-odd-even Hecke operators $R$ and $R'$ respectively, are equivalent as braided abelian categories if and only if $R$ and $R'$ are defined for the same $q$ and are of the same super rank. 
\end{edl}

As in the case of even Hecke operators, we have
\begin{cor}\label{cor52} The super rank of an odd-even Hecke operator and the dimension of $V$, on which $R$ operates, satisfies the inequality $m+n\leq \dim(V)$, and if the equality takes place then $P_\Lambda(t)=(1+t)^m(1-t)^{-n}$, hence $P_S(t)=(1+t)^n(1-t)^{-m}$, and by Equation \rref{eq23}, $P_\Lambda(t)=(1+t)^{2mn}(1-t)^{-m^2-n^2}$.
\section{The Poincar\'e series of $E$}\label{secE}\end{cor}

We now give a relation between the Poincar\'e series of $E$ and the ones of $\Lambda$ and $S$. The formula given here has been proved for the case $q=1$ by Lyubashenko \cite{lyu} and $q$ transcendent by the author \cite{ph97}. Here, using the general theory of symmetric functions, we recover this relation for any $q$ not root of unity. According to \rref{eq10}, 
\bbas e_n=\dim E_n=\sum_{\lambda\part n}(m_\lambda)^2.\eeas
On the other hand, we have the following identity of Schur functions, for two countable set of variables $x=(x_1,x_2,\ldots), y=(y_1,y_2,\ldots)$,
\bbas \prod_{i,j}(1-x_iy_j)^{-1}=\sum_{\lambda\in {\cal P}}\bfm_\lambda(x)\bfm_\lambda(y).\eeas
Setting $y_i=tx_i$, we have
\begin{equation}\label{eq22} \prod_{i,j}(1-x_ix_jt)^{-1}=\sum_{\lambda \in {\cal P}}\bfm_\lambda(x)^2t^{|\lambda|}.\end{equation}
The left-hand side of this equation can be considered as a generating function for complete symmetric functions in variables $z_{ij}=x_ix_j$. Setting $\bfm_\lambda=m_\lambda$ then the right-hand side of \rref{eq22} is equal to $P_E(t)$, by \rref{eq10}. Thus we derive from \rref{eq22}
\begin{equation}\label{eq23}P_E(t)=P_S(t)\star P_S(t),\end{equation}
where $\star$ denotes the multiplication on the $\lambda$-ring $\BBb{K}_0[[t]]$ of formal power series with first coefficient equal to 1 \cite{knutson}. If $P(t)=\prod(1+x_it), Q(t)=\prod(1+y_it)$ then $P\star Q(t):=\prod(1+x_iy_jt)$. We shall refer to $\star$ as the $\lambda$-product. Recall that the $\BBb{K}_0[[t]]$ is a ring with the above $\lambda$-product and the $\lambda$-addition is the usual product of power series [loc.cit.].

Let $R_1,R_2$ be Hecke operators on $V_1,V_2$, their Hecke sum $R$ acts on $V_1\oplus V_2$ as follows, $R|V_i\ot V_i=R_i$, $R(v_1\ot v_2)=qv_2\ot v_1, R(v_2\ot v_1)=v_1\ot v_2-(q-1)v_2\ot v_1$, for $v_i\in V_i$. Then $R$ is a Hecke operator, too. Let $\Lambda_R, S_R$ be the associated quadratic algebras, then it is not difficult to check that
\bbas P_{\Lambda_R}(t)=P_{\Lambda_{R_1}}(t)P_{\Lambda_{R_1}}(t),\\
P_{S_R}(t)=P_{S_{R_1}}(t)P_{S_{R_1}}(t).
\eeas
Thus, the Hecke sum gives rise to the $\lambda$-addition on $\BBb{K}_0[[t]]$.

Let $P(t):=\frac{d}{dt}{\rm ln}(P_S(t))=P'_S(t)P_S(t)^{-1}$. $P(t)$ is the generating function for the power sums (\cite{mcd1,mcd2}, I.3). The operator $\frac{d}{dt}$ maps $\BBb{K}_0[[t]]$ bijectively to $\BBb{K}[[t]]$ -- the ring of power series, the inverse operator is $\int_0^t$. The $\lambda$-addition is mapped to the ordinary addition while the $\lambda$-multiplication is mapped to the component-product (denoted by $*$). Indeed, let $\bfp_r(x)$ (resp. $\bfp_r(y)$) be the $r^{\rm th}$ power sum of $x=(x_1,x_2,\ldots)$ (resp. $y=(y_1,y_2,\ldots)$), then $\bfp_r(xy)=\bfp_r(x)\bfp_r(y),$ where $xy:=\{x_iy_j|1\leq i,j\}$. Thus we have (cf. \cite{ph97}),
\begin{equation}P_E(t)=\exp\int_0^tP(u)^{*2}.\end{equation}

{\bf Acknowledgement}
\ictpthanks\ 
It is a pleasure to thank Professors I.G. Macdonald and J. Stembridge  for their help.

%\bibliography{bible}

\end{document}